\documentclass[11pt]{article}

\usepackage{amssymb}
\usepackage{latexsym}
\usepackage[mathscr]{euscript}
\usepackage{graphicx}
\usepackage{float}
\usepackage{subfigure}

\usepackage{amscd}
\usepackage{amsmath}
\usepackage{amssymb}
\usepackage{fullpage}

\usepackage{amssymb}
\usepackage{latexsym}
\usepackage[mathscr]{euscript}
\usepackage{pb-diagram}
\usepackage{latexsym}
\usepackage{amsfonts}
\usepackage{epsfig}

\newcommand\bb[1]{\Bbb{#1}}

\newcommand\pr{^{\prime}}

\makeatletter
\def\bbordermatrix#1{\begingroup \m@th
  \@tempdima 4.75\p@
  \setbox\z@\vbox{%
    \def\cr{\crcr\noalign{\kern2\p@\global\let\cr\endline}}%
    \ialign{$##$\hfil\kern2\p@\kern\@tempdima&\thinspace\hfil$##$\hfil
      &&\quad\hfil$##$\hfil\crcr
      \omit\strut\hfil\crcr\noalign{\kern-\baselineskip}%
      #1\crcr\omit\strut\cr}}%
  \setbox\tw@\vbox{\unvcopy\z@\global\setbox\@ne\lastbox}%
  \setbox\tw@\hbox{\unhbox\@ne\unskip\global\setbox\@ne\lastbox}%
  \setbox\tw@\hbox{$\kern\wd\@ne\kern-\@tempdima\left[\kern-\wd\@ne
    \global\setbox\@ne\vbox{\box\@ne\kern2\p@}%
    \vcenter{\kern-\ht\@ne\unvbox\z@\kern-\baselineskip}\,\right]$}%
  \null\;\vbox{\kern\ht\@ne\box\tw@}\endgroup}
\makeatother

\newcommand\commentout[1]{\marginpar{\tiny $\backslash$commentout}}

\newcommand\qed{\hfill$\square$}
\def\column#1#2{\mathrel{\mathop{#1}\limits_{#2}}}

\def\compcirc {\mbox{\hspace{.05cm}}\raisebox{.04cm}{\tiny  {$\circ$ }}}

\newtheorem{Lemma}{Lemma}[section]
\newtheorem{Theorem}{Theorem}[section]
\newtheorem{Proposition}{Proposition}[section]
\newtheorem{Definition}{Definition}[section]

\newtheorem{Example}{Example}[section]
\newtheorem{Remark}{Remark}[section]

\newenvironment{Proof}{\par\noindent\textbf{Proof:}}
{\qed}

\usepackage{graphics}
\title{Approximating Absolute Galois Groups}

\author{Gunnar Carlsson \\ Department of Mathematics, Stanford University  \\Stanford, California 94305 \and Roy Joshua \\Department of Mathematics, The Ohio State University \\Columbus, Ohio 43210}
\begin{document}
\maketitle
{\bf Abstract:} In this paper we identify a class of profinite groups (totally torsion free groups)  that  includes all separable Galois groups of fields containing an algebraically closed subfield, and demonstrate that it can be realized as an inverse limit of torsion free virtually finitely generated abelian (tfvfga)  profinite groups.  We show by examples that  the condition is quite restrictive.  In particular, semidirect products of torsion free abelian groups are rarely totally torsion free.  The result is of importance for $K$-theoretic applications, since descent problems for tfvfga  groups are relatively manageable.  
\section{Introduction}
It is well understood that the structure of absolute Galois groups is quite restricted, and that general profinite groups cannot be absolute Galois groups. Here are some examples of the restrictions that are known or conjectured.  
\begin{enumerate}
\item{The {\em Artin-Schreier theorem} asserts that the only finite groups that can occur as subgroups of absolute Galois groups are the trivial group and the cyclic group of order two.  }
\item{ F. Bogomolov has conjectured in \cite{bogomolov} that in the case where the base field contains an algebraically closed field,  the $p$-Sylow subgroup of the commutator subgroup is free.  This is obviously a very restrictive condition.  }
\item{ The  Bloch-Kato conjecture \cite{blochkato} implies in particular  that when the roots of unity are in the base field, then the cohomology ring of the absolute Galois group is generated in degree one and is defined by relations in degree two.  This is also clearly a very restrictive condition.}
\end{enumerate} 
In this paper, we study another condition which applies to all absolute and separable Galois groups.  Specifically, we prove three  results.  

\begin{enumerate}
\item{We identify a new group theoretic  condition, {\em total torsion freeness} (see Definition \ref{ttfdef}), and prove that it holds for absolute Galois groups of fields which contain all roots of unity. The condition can be interpreted both in the case of discrete groups and profinite groups.  Total torsion  freeness condition implies the usual notion of torsion freeness.    }
\item{We prove that all groups satisfying total torsion freeness can be approximated in an appropriate sense by groups in a very restricted family which we call $\frak{B}$. In the case of profinite groups, this means that the group can be described as an inverse limit of groups in $\frak{B}$.    The family $\frak{B}$ consists of groups which are (a) virtually finitely generated abelian (virtually topologically finitely generated abelian in the profinite case) and (b) torsion free.  See Definition \ref{defvirt} for he definition.  }
\item{We  prove that every group in the family $\frak{B}$ embeds as a torsion free subgroup of the group $\Sigma _n \ltimes \bb{Z}^n $ (or $\Sigma _n \ltimes \hat{\bb{Z}}^n$ in the profinite case) for some $n$.  We also show that it follows that every group $\Gamma$ in $\frak{B}$ acts on the ring 
$$ \bb{E}_n = \bigcup _{s} k[t_1^{\pm \frac{1}{s}}, \ldots ,  t_n^{\pm \frac{1}{s}}]$$
for some $n$,  in such way that the extension $\bb{E}_n^{\Gamma} \subseteq \bb{E}_n$ is a Galois ring extension in the sense of \cite{magid}.  This last statement will be extremely useful in our work on the algebraic $K$-theory of fields.  }
\end{enumerate} 

Example of groups satisfying the total torsion freeness condition include free groups, free profinite groups, free abelian groups, free abelian profinite groups, fundamental groups of orientable two manifolds (not necessarily compact), and free products of totally torsion free groups.  On the other hand, groups that are torsion free in the usual sense but not totally torsion free include fundamental groups of non-orientable surfaces, and the groups of upper triangular integral $n \times n$ matrices with ones along the diagonal.  It is also easy to check  that most semidirect products of the form $\bb{Z} \ltimes \bb{Z}^n$ fail to be totally torsion free but are of course torsion free.  

This work is motivated by our work on the descent problem in the algebraic $K$-theory of fields.  We will  not go into detail on this work here, but will note that a key component is understanding the structure of an analogue of the geometric classifying space construction of \cite{morelvoevodsky} for profinite groups, and the approximation theorem in this paper  permits the construction of very useful and explicit models of such classifying spaces.   We believe that the results are also of independent 
interest, hence the present  paper. 

The first author wishes to express his thanks to Brian Conrad for a number of stimulating conversations, in particular concerning the proof of Theorem \ref{property}.  
\section{Total torsion freeness} In this section, we define total torsion freeness and prove that separable Galois groups of fields have this property. 
\begin{Definition} \label{ttfdef} A Hausdorff topological  group $G$  is said to be {\em totally torsion free} if the abelianization of any closed subgroup of $G$ is torsion free. In this case, the abelianization of a subgroup $K$ means the quotient of $K$ by the closure of its commutator subgroup.   A topological  group is said to be {\em weakly totally torsion free} if every closed subgroup of finite index  has torsion free abelianization. 
\end{Definition}
\begin{Remark} {\em The only groups with non-discrete  topology we study will be profinite groups.  The general statement is made only since it allows uniform treatment of the discrete and profinite cases. }
\end{Remark}  
\begin{Remark}{\em  Note that total torsion freeness implies torsion freeness, since a group is torsion free if and only if every cyclic subgroup is torsion free and, since it is abelian, its abelianization is torsion free.  In the profinite case, we have the same situation, since topologically cyclic subgroups must be torsion free since they are abelian.  }
\end{Remark}
\begin{Example} Free groups and free abelian groups are totally torsion free.  Similarly for free and free abelian profinite and pro-$l$ groups. See Proposition \ref{freeprod} for a proof. 

\end{Example}

\begin{Example} The integral Heisenberg group is torsion free, but not totally torsion free.  The group consists of upper triangular integral matrices with diagonal elements all equal to one.  It is easy to check that the subgroup consisting of all elements of the form 
$$
\begin{bmatrix}
1  &  2m    &  n  \\
0  &  1  &  l  \\
0  &  0  &  1 
\end{bmatrix}
$$

has abelianization isomorphic to $\Bbb{Z} \times \Bbb{Z} \times \Bbb{Z}/2\Bbb{Z}$, and that therefore the group is not totally torsion free.   
\end{Example}
\begin{Example} Fundamental groups of non-orientable surfaces are torsion free but not totally torsion free.  They are torsion free since they act freely on a finite dimensional contractible manifold, but not totally torsion free since their one-dimensional integral homology has torsion.  
\end{Example} 
One elementary result about this notion is the following.  
\begin{Proposition} \label{freeprod} Free products of totally torsion free groups are totally torsion free.   Free products of weakly totally torsion free profinite (pro-$p$) groups are weakly totally torsion free.  
\end{Proposition} 
\begin{Proof} The first statement follows immediately from the Kurosh theorem. The second statement follows from the analogous  profinite and pro-$p$ statements in \cite{binz} and \cite{gildenhuys}.  
\end{Proof}

We also have the following.  
\begin{Proposition} The fundamental groups of orientable two dimensional manifolds $M$ without boundary are totally torsion free.  
\end{Proposition}
\begin{Proof} It suffices to prove that $H_1(M, \bb{Z}) \cong \pi _1(M)^{ab}$ is torsion free, since any subgroup of $\pi _1 (M)$ is also the fundamental group of an orientable 2-manifold without boundary.  But this result is Corollary 7.13 of \cite{bredon}.  
\end{Proof} 

We will also need to record results concerning Pontrjagin duality.  Recall that for a compact topological abelian group $A$, we define  the Pontrjagin dual to $A$, denoted $\hat{A}$,  to be $Hom^c(A,S^1)$, where $S^1$ denotes the circle group, and the superscript ``c" denotes continuous homomorphisms.  

\begin{Proposition} \label{pontrjagin} The construction $A \rightarrow \hat{A}$ satisfies the following properties. 
\begin{enumerate}
\item{\label{intro}The $\hat{}$-construction defines an equivalence of categories from the category of compact topological abelian groups to the opposite of the category of discrete abelian groups. The $\hat{}$-construction is its own inverse.  }
\item{ \label{discrete} For a profinite group $G$, $\hat{G}$ is  isomorphic to $Hom^c(G,\mu _{\infty})$, where $\mu _{\infty} \subseteq S^1$ is the group of all roots of unity, isomorphic to $\bb{Q}/\bb{Z}$.  If $G$ is a $p$-profinite group, then $\mu _{\infty}$ can be replaced by $\mu _{p^{\infty}}$, the group of all $p$-power roots of unity, isomorphic to $\bb{Z}[\frac{1}{p}]/\bb{Z}$.  }
\item{The functor $A \rightarrow \hat{A}$ is exact.  }
\item{ \label{criterion} For $G$ a profinite abelian group, $G$ is torsion free if and only if $\hat{G}$ is divisible.  Similarly for ``$p$-torsion free" and ``$p$-divisible".  }
\end{enumerate} 
 \end{Proposition}
 \begin{Proof} Statement (1) is one version of the statement of the Pontrjagin duality theorem, (2) is an immediate consequence, and (3) follows immediately from (1). It remains to prove (4).  To prove (4), we note that $G$ is torsion free if and only if the sequence $0 \rightarrow G \stackrel{\times n}{\longrightarrow} G$ is exact.  The exactness proves that this occurs if and only if $\hat{G} \stackrel{\times n}{\longrightarrow} \hat{G} \rightarrow 0$ is exact, so $\times n$ is surjective.   This is  the result.  

 \end{Proof}

We now have the main result of this section.

\begin{Theorem} \label{property} Let $F$ be any field containing all roots of unity.  Then the absolute  Galois group $G_F$ of $F$ is totally torsion free.  
\end{Theorem} 

\begin{Remark} {\em Class field theory shows, for example,  that one cannot expect this result to hold for absolute Galois groups of number fields, so that some condition on the field is necessary.  } 
\end{Remark}
\begin{Proof}  Consider any closed subgroup $K \subseteq G_F$, and its corresponding extension $F_K$ of $F$.  We are interested in the Galois group $\frak{A} = G_{F_K}^{ab}$  of the extension $F_K^{ab} $ over $F_K$, and want to prove that it is torsion free.  We have that 
$$\frak{A} \cong \mbox{\hspace{.1mm}}\column{lim}{\column{\leftarrow}{n}} \frak{A}/n\frak{A}
$$
where $n$ varies over the partially ordered set of integers with ordering given by  $n \leq n\pr$ if and only if $n|n\pr$.  For each prime, we let $\frak{A}_p$ denote the $p$-Sylow subgroup of $\frak{A}$.  Of course, we have $\frak{A} \cong \prod _p \frak{A}_p$, and it will suffice to prove that $\frak{A}_p$ is torsion free for each $p$.  When $p \neq \mbox{char}(F)$, we can apply Kummer theory as follows.    Kummer theory asserts that 
$$Hom^c(\frak{A}/{p^k}\frak{A}, \Bbb{Z}/p^k\Bbb{Z}) \cong F_k^*/p^kF^*_K \cong F^* \otimes C_{p^k}$$ where 
 $C_{p^n}$ denotes the cyclic group of order $p^n$.  We can pass to the direct limit on both sides of this isomorphism to get an isomorphism   
$$  Hom^c(\frak{A}_p, C_{p^{\infty}} ) \rightarrow F^*_K \otimes C_{p^{\infty}}
$$ 
where $C_{p^{\infty}}$ denotes the group  $\Bbb{Z}[\frac{1}{p}]/\Bbb{Z}$.  We note that the group $ Hom^c(\frak{A}_p, C_{p^{\infty}} )$ can be interpreted as the Pontrjagin dual  $\hat{\frak{A}_p}$, by part \ref{discrete} of Theorem \ref{pontrjagin}. 
The discrete group $F_K^* \otimes   C_{p^{\infty}}$ is clearly $p$-divisible, so by part \ref{criterion}  of Proposition \ref{pontrjagin}, $\frak{A}_p$  is torsion free. 

On the other hand, if $p = \mbox{char}(F)$, we use Witt vectors instead, as in \cite{lorenz}, Ch. 26. We recall that for any algebra $A$ over the finite field $\Bbb{F}_p$, we can construct the ring of {\em Witt vectors} $W(A)$.  It has a number of useful properties.  
\begin{enumerate} \item{As a set, $W(A)$ is the infinite product $\prod _{n=0}^{\infty} A$}.  
\item{The shift operator 
$$ V(a_0,a_1, \ldots ) = (0,a_0,a_1, \ldots )
$$
is a homomorphism of abelian groups.}
\item{The {\em Frobenius operator} $F$  defined by $F(a_0, a_1, \ldots )  = (a_0^p, a_1^p, \ldots)$ commutes with $V$, and the composite $VF = FV$ is multiplication by $p$ in the group structure on $W(A)$.  }
\item{The subset $I_n(A) \subseteq W(A)$ given by $I_n(A) = V^nW(A)$ is an ideal in $W(A)$, and we denote the quotient algebra $W(A)/I_n(A)$ by $W_n(A)$.  Multiplication by $V$ gives an operator $W_n(A) \rightarrow W_{n+1}(A)$.  }
\end{enumerate}
We now use these constructions to study the dual to $\frak{A}_p$, following F. Lorenz \cite{lorenz}. Let $\pi:W(F) \rightarrow W(F)$ denote the operator $F - id$.  Because $F$ and $V$ commute, it is clear that $V: W_n(F) \rightarrow W_{n+1}(F)$ induces  a  homomorphism 
$$  \tilde{F}: W_n(F)/\pi (W_n(F)) \rightarrow W_{n+1}(F)/\pi (W_{n+1}(F)) 
$$
and what is proved in \cite{lorenz}, Ch. 26, p. 108,  is that there is a perfect duality between the groups 
$$W_n(F_K)/\pi (W_n(F_K)) \times \frak{A}_p/p^n \frak{A}_p \rightarrow C_{p^{\infty}}
$$
and consequently an isomorphism $\column{colim}{\rightarrow}W_n(F_K)/\pi (W_n(F_K))  \rightarrow Hom^c(\frak{A}_p, C_{p^{\infty}})$, after passage to colimits over $n$.    We claim that the group $\column{colim}{\rightarrow}W_n(F_K)/\pi (W_n(F_K)) $ is $p$-divisible, from which the torsion freeness of $\frak{A}_p$ would follow as in the case $p \neq \mbox{char}(F)$ above.  To see this, one must only observe that $FV$ also induces a map 
$$  \tilde{FV}: W_n(F)/\pi (W_n(F)) \rightarrow W_{n+1}(F)/\pi (W_{n+1}(F)) 
$$
and that because $x = F(x) $ in $W_n(F)/\pi (W_n(F)) $, the two maps $\tilde{F}$ and $\tilde{FV}$ are equal.  Since $FV$ is multiplication by $p$, the result follows.  This gives the required result.  
\end{Proof}
\section{The approximation theorem} This section will define a class of groups which can be used to approximate any totally torsion free profinite group, in the sense that the group can be described as an inverse limit of groups in that family.  
\begin{Definition} \label{defvirt} A profinite  group $G$ is said to be {\em virtually finitely generated abelian} if it contains a closed topologically finitely generated abelian subgroup of finite index.  We will denote the family of all virtually finitely generated abelian profinite groups by ${\cal{V}}^{ab}$.  The family of torsion free groups in ${\cal V}^{ab}$ will be denoted by $\frak{B}$. In the discrete situation, we will denote by $\frak{B}^{0}$ the family of all torsion free groups that contain a finitely generated abelian subgroup of finite index.  
\end{Definition} 
\begin{Example} {\em Let $A$ and $B$ be finitely generated torsion free abelian groups, and let $\varphi: A \longrightarrow Aut(B)$ be a homomorphism.  Suppose further that the image of $\varphi$ is finite.  Then the semidirect product $A \ltimes _{\varphi} B$ is torsion free and virtually finitely generated abelian.  The torsion free property holds for any semidirect product of one free abelian group with another, and the fact that it is virtually finitely generated abelian follows from the observation that $Ker(\varphi) \times B \subseteq A \ltimes _{\varphi}B$ is a normal subgroup of finite index, and therefore $A \ltimes _{\varphi}B$ belongs to $\frak{B}^0$. The profinite completion of any such group gives a profinite group in $\frak{B}$.   }
\end{Example} 
\begin{Example} {\em Let $G$ be the fundamental group of any compact  flat Riemannian manifold.  Then by the results of \cite{auslander} $G$ is in the family $\frak{B}^0$.  The profinite completion of $G$ gives a profinite group in $\frak{B}$. } 
\end{Example} 
\begin{Remark}{\em Note that these examples are not {\em totally} torsion free.  Although groups in $\frak{B}$  will be used to approximate totally torsion free groups, they are not generally totally torsion free themselves.  }
\end{Remark} 
It is easy to see that if a profinite group is virtually abelian, then it  admits a closed {\em normal} abelian subgroup of finite index.  This follows since if $K$ is any finite index abelian subgroup of a profinite group, then the intersection of the (finite) collection of conjugates of $K$ is the required abelian normal subgroup.  Our initial goal is to prove that given a totally torsion free profinite group $\Gamma$, and a continuous surjective homomorphism 
$\pi: \Gamma \rightarrow G$, where $G$ is finite, then there is a profinite group $\overline{G} \in \frak{B}$ and  surjections $\overline{\pi}: \Gamma \rightarrow \overline{G}$ and $\sigma: \overline{G} \rightarrow G$, with the kernel of $\sigma $ abelian and finitely generated,  so that $\sigma \compcirc \overline{\pi} = \pi$.  So, we fix a finite quotient group of $G$ of $\Gamma$, and define an {\em approximation system} for $\Gamma \rightarrow G$ to be a sequence of  continuous surjective homomorphisms 
$$  \Gamma \longrightarrow \hat{G} \stackrel{\sigma}{\longrightarrow} G 
$$
where the composite is the projection from $\Gamma$ to $G$, and where the kernel of $\sigma $ is torsion free, finitely generated,  and abelian.

Let $p$ be a prime,  let $\alpha:\Gamma \rightarrow \hat{G} \rightarrow $G be an approximation system for $\Gamma \rightarrow G$, and let $g \in G$ be of order $p$.  We say that   $\alpha$ is {\em $p$-torsion free over $g$} if  $\sigma ^{-1}(g) \subseteq \hat{G}$ contains no $p$-torsion elements.  

\begin{Proposition}  Let $\pi:\Gamma \rightarrow G$ be a continuous homomorphism of profinite groups, where $G$ is  finite, and $\Gamma $ is weakly totally torsion free.  Then for any prime $p$ and element $g$ of order $p$  in $G$, there is an  approximation system
$$  \Gamma \longrightarrow \hat{G}_g \stackrel{\sigma _g}{\longrightarrow } G
$$
which is $p$-torsion free over $g$.  
\end{Proposition}
\begin{Proof} Let $g$ be an element of order $p$ in $G$.   The subgroup $L = \pi ^{-1}  \langle g \rangle \subseteq \Gamma$ is  closed, and therefore has torsion free abelianization $L^{ab}$.  The restriction of $\pi $ to $L$ gives a surjective homomorphism from $L$ to $\langle g \rangle \cong \Bbb{Z}/p\Bbb{Z}$, which we also denote by $\pi$.  Let $L^{ab}_p$ denote the $p$-Sylow subgroup of $L^{ab}$, which is naturally both a subgroup and a quotient of $L^{ab}$, and let $r_p$ denote the projection $ L^{ab} \rightarrow L^{ab}_p$.  Since $\pi$ vanishes on all $q$-Sylow subgroups (for $q \neq p$) of $L^{ab}$, it follows that $\pi$ naturally factors through a homomorphism $\pi _p: L^{ab}_p \rightarrow \Bbb{Z}/p\Bbb{Z}$.  The Pontrjagin dual $\frak{D} = Hom ^c (L^{ab} _p, \mu _{p^{\infty}})$, where $\mu_{p^{\infty}} \cong \Bbb{Z}[\frac{1}{p}]/\Bbb{Z}$, is a divisible group because $ L^{ab} _p$ is torsion free.  The homomorphism $\pi _p$ is  identified with a non-zero element of order $p$  in  $\frak{D}$.  Because of the divisibility of $\frak{D}$, there is a sequence of elements $\varphi _i \in \frak{D}$ with $\varphi _1 = \pi _p$, and $p \cdot \varphi _i = \varphi _{i-1}$, or equivalently, a homomorphism $\mu _{p^{\infty}}$ into $\frak{D}$ extending $\varphi _1$.  Since $\mu _{p^{\infty}} $ is the Pontrjagin dual to $\Bbb{Z}_p$, we obtain a surjective homomorphism $\Pi_p: L^{ab}_p \rightarrow \Bbb{Z}_p$, which projects to $\pi _p$ in $\Bbb{Z}/p\Bbb{Z}$.  Let $\overline{\Pi}_p$ denote the composite 
$$  L \rightarrow L^{ab}_p \stackrel{\Pi _p}{\longrightarrow} \bb{Z}_p
$$

 The kernel of $\overline{\Pi}_p$ (denote it by $\frak{L}$) is a closed subgroup of $L$, and the quotient 
$L /\frak{L}$ is isomorphic to $\Bbb{Z}_p$.  Of course, $\frak{L}$ is also a closed subgroup of $ K \subseteq \Gamma$, where $K$ is the kernel of $\pi$, but it is not necessarily normal in $\Gamma$.  However, it has only finitely many conjugates, as  it is normalized by $L$, which has finite index in $\Gamma$.  Let $\{\frak{L}_1, \ldots , \frak{L}_n \}$ be the set of conjugates of $\frak{L}$, and let $\overline{\frak{L}}$ denote the intersection 
$$ \bigcap _{i=1}^n \frak{L}_i
$$
The closed subgroup $\overline{\frak{L}}$ is normal in $\Gamma$, and we consider the quotient 
$\Gamma /\overline{\frak{L}}$.  We have the short exact sequence of groups 
$$ \{ e \} \rightarrow K/\overline{\frak{L}} \rightarrow \Gamma /\overline{\frak{L}} \rightarrow G \rightarrow \{ e \}
$$
Since $K/\overline{\frak{L}} \subseteq \prod _i K/\frak{L}_i$ we have that $K /\overline{\frak{L}}$ is a torsion free profinite group, and that it is topologically finitely generated.  Finally, we need to show that 
$\Gamma /\overline{\frak{L}} $ is $p$-torsion free over $g$.  But it is clear that any lift of $g$ to $L$ projects to the image of $g$ under the homomorphism $\pi _p$, and is consequently a topological generator for $\bb{Z}_p$.  In particular, it cannot be a torsion element.    \end{Proof}

\begin{Definition} Let $\alpha _1: \Gamma \rightarrow \hat{G}_1 \stackrel{\sigma _1}{\longrightarrow} G$ and $\alpha _2: \Gamma \rightarrow \hat{G}_2 \stackrel{\sigma _2}{\longrightarrow} G$
be two approximation systems for $\Gamma \rightarrow G$.  Let $\hat{G}_{12}$ denote  the image of $\Gamma$ in $\hat{G}_1 \times
 \hat{G}_2$. 
  It is clear that the composites $\hat{G}_{12} \stackrel{\pi _1}{\longrightarrow} \hat{G}_1\stackrel{\sigma_1}{ \longrightarrow }G$ and $\hat{G}_{12} \stackrel{\pi _2}
{\longrightarrow} \hat{G}_2 \stackrel{\sigma _2}{\longrightarrow }G$
  are equal, so we obtain a homomorphism $\sigma _{12} :\hat{G}_{12} \rightarrow G$.  We will define the {\em fiber product} of $\alpha _1$ and $\alpha _2$ to be the approximation system 
$$  \Gamma \longrightarrow \hat{G}_{12} \stackrel{\sigma _{12}}{\longrightarrow } G
$$
and denote it by $\alpha_1 \column{\times }{G} \alpha _2$. 
\end{Definition}

The following proposition will allow us to construct an approximation system $\Gamma \rightarrow \hat{G} \rightarrow G$ for which $\hat{G}$ is finitely generated and  torsion free.  

 \begin{Proposition}Let $\alpha _1:\Gamma \longrightarrow \hat{G}_1 \stackrel{\sigma _1}{\longrightarrow } G$ and $\alpha _2:\Gamma \longrightarrow \hat{G}_2 \stackrel{\sigma _2}{\longrightarrow } G$ be approximation systems for $\Gamma \rightarrow G$, where $\Gamma $ is profinite and $G$ is finite.  Suppose that $S_1$ and $S_2$ are finite collections of pairs $(p,g)$, where $p$ is a prime and $g \in G$ is an element of order $p$. Suppose further that the approximation system $\alpha _i$ is such that for any $(p,g) \in S_i$, $\alpha _i$ is $p$-torsion free over $g$.    Then $\alpha _1 \column{\times}{G} \alpha _2$ is $p$ torsion free over $g$ for any $(p,g)  \in S_1 \cup S_2$.  
\end{Proposition}
\begin{Proof} Clear from the definitions.  
\end{Proof}

\begin{Theorem} \label{main}Let $\Gamma $ be a weakly totally torsion free profinite group, and suppose we are given a continuous homomorphism $\Gamma \rightarrow G$.  Then there is a an approximation system 
$$ \Gamma \longrightarrow \hat{G} \rightarrow G
$$
with $\hat{G}$ torsion free.  
\end{Theorem}
\begin{Proof} Let $\{ (p_1,g_i), \ldots , (p_m,g_m) \}$ be a complete list of all pairs $(p,g)$ so that $p$ is a prime and $g$ is a $p$-torsion element of $G$.  For each $i$, construct an approximation system 
$$  \alpha _i :\Gamma \longrightarrow \hat{G}_ i\rightarrow G
$$
so that $\hat{G}_i$ is $p$-torsion free over $g_i$.  Form the iterated fiber product 
$$ \alpha _1 \column{\times}{G} \cdots \column{\times}{G} \alpha_m = \{\Gamma \rightarrow \hat{G} \rightarrow G \}
$$
It is $p$-torsion free for every element $g$ of order $p$ in $G$.  Suppose that there is a torsion element $g$ in $\hat{G}$.  Then by choosing one of the primes $p$ dividing the order of $g$, and taking appropriate powers, we can generate an element $g^{\prime}$ of order $p$.  Since the kernel of the projection $\hat{G} \rightarrow G$ is a torsion free group, $g^{\prime}$  must have a non-trivial projection in $G$.   But from the construction, its powers will all be non-identity in the fiber product factor corresponding to $(p,g^{\prime})$, which gives the result.  
\end{Proof}

This implies the following result concerning weakly totally torsion free profinite groups. 

\begin{Theorem} Let $\Gamma$ denote a weakly totally torsion free profinite group.  Let $\frak{V}$ denote the set of all closed normal subgroups $K$ of $\Gamma$ for which $\Gamma /K$ is in $\frak{B}$.  Then we have the natural homomorphism
$$
\varphi:\Gamma \longrightarrow \prod _{K \in \frak{V}} \Gamma/K
$$
The map $\varphi$ is an inclusion and the group $\Gamma$ is isomorphic as a topological group to the image of $\varphi$. Equivalently, $\Gamma$ is isomorphic to the inverse limit of the partially ordered set of quotients in $\frak{V}$.  
\end{Theorem} 
\begin{Proof} The injectivity of $\varphi$ is an immediate  consequence of Theorem \ref{main}, and the statement about being homeomorphic to its image is a standard property of compact Hausdorff spaces.  The inverse limit statement is similarly elementary.  
\end{Proof} 
\begin{Remark}{\em Note that the approximating groups $\Gamma /K$ are only required to be torsion free, they are likely not totally torsion free themselves.   }
\end{Remark} 

There is an analogous statement for discrete groups, which can be proved in an entirely analogous manner.  Since the final conclusion is weaker, we state it without proof. 

\begin{Theorem} Let $\Gamma$ we a weakly totally torsion free discrete group, and suppose further that $\Gamma$ is residually finite.  Then for any element $\gamma \in \Gamma$, there is a homomorphism $f:\Gamma \rightarrow G$, where $G \in \frak{B}^{0}$, so that $f(\gamma) \neq e$.  
\end{Theorem}

\section{Galois theoretic properties}
We will begin by showing that every group in $\frak{B}$ can be embedded as a closed subgroup of the semidirect product  $\Sigma _n \ltimes \hat{\bb{Z}}^n$ for some $n$. As a set, $\Sigma _n \ltimes \hat{\bb{Z}}^n$ is isomorphic to the product $\Sigma _n \times \hat{\bb{Z}}^n$, and we topologize it as such, with the topology on $\Sigma _n$ being the discrete topology.  It is then clearly Hausdorff, compact, and totally disconnected, and it is easy to check that the multiplication remains continuous, from which it follows that it is a profinite group.  We first recall the definition of the wreath product. 

\begin{Definition} \label{wreathdefinition}   Let $G$ denote a discrete group, and $K$ any group.  Let $K^G$ denote the set of all functions $f :G \rightarrow K$, made into a group by equipping it with pointwise multiplication.  Also, we equip  $K^G$ with the left  $G$-action by automorphisms  $(g \cdot f) (g^{\prime}) = f(g^{\prime}g^{-1})$.  We define the {\em wreath product} $G \wr K$ to be the semidirect product $G \ltimes K^G$ with the given action.  We note that in the situation where $K$ is a topological group (but $G$ remains discrete), this construction still makes sense, and $G\wr K$ is in a natural way a topological group.  
\end{Definition} 
We  consider the situation of a topological group $G$ equipped with a continuous surjective homomorphism $\pi: G \rightarrow Q$, where $Q$ is equipped with the discrete topology, and where the kernel of $\pi$ is abelian.  
\begin{Lemma} \label{wreathembed} Suppose that $G$ is as above, and  is Hausdorff.  Let $K$ denote the  kernel of $\pi$, an abelian Hausdorff topological group.  Then there is a closed  embedding $i: G \rightarrow Q \wr K$ over $Q$, in the sense that the diagram 
$$ 
\begin{diagram} \node{G} \arrow{e,t}{i} \arrow{se,t} {\pi} \node {Q \wr K} \arrow{s} \\
\node{}\node{Q}
\end{diagram} 
$$
commutes, so the image of $i$ is a closed subgroup of $Q \wr K$.  
\end{Lemma}
\begin{Proof} We first consider the case where $G$ is discrete.  
The standard classification of extensions of groups with abelian kernel shows that group structures on the set $Q \times K$ for which the projection $Q \times K \rightarrow Q$ and the inclusion $K \hookrightarrow Q \times K$ given by $k \rightarrow (e,k)$  are homomorphisms  are in one to one correspondence with $2$-cocycles on $G$ with values in the $G$-module $K$, and that cocycles $c_1$ and $c_2$ determine isomorphic groups if $c_1 - c_2$ is   a coboundary.  For a $2$-cocycle $c$, we let $G(c)$ denote the corresponding group.  It is immediate that given a group $Q$ and a homomorphism of $G$-modules $f:K_1 \rightarrow K_2$, there is a naturally associated homomorphism $G(c) \rightarrow G(f\compcirc c)$, which respects the projection to $Q$, and for which the restriction of the homomorphism to the kernels is the homomorphism $f$.      Thus our exact sequence $1 \rightarrow K \rightarrow G \rightarrow  Q \rightarrow 1$ is associated to a  $2$-cocycle $c$  on $Q$ with values in $K$.  We consider the induced $Q$-module $K^Q$ as in Definition \ref{wreathdefinition} above.   Shapiro's lemma asserts that $H^2(Q,K^Q) $ vanishes, and we see that for any $2$-cocycle $c$  on $Q$ with values in $K^Q$, the group $G(c)$ is isomorphic to the wreath product $Q  \wr K$.   Let $c$ denote a cocycle defining the extension $1 \rightarrow K \rightarrow G \rightarrow Q \rightarrow 1$.  There is a natural homomorphism   $i_K:K \rightarrow K^Q$ of $Q$-modules, which sends $k \in K$ to the constant function on $Q$ with value $k$.  We therefore obtain a homomorphism $G(c) \rightarrow G(i_K\cdot c)$ of groups, suitably compatible with the projections to $Q$ and the inclusions of the kernels, and from the above discussion   $G(i_K \cdot c)$ is isomorphic to the wreath product $Q \wr K$.  In the non-discrete case, the identical method works provided one verifies continuity of the  homomorphism $i_K$, which is immediate.  
\end{Proof}

Suppose now that we have a group $G$ in the family $\frak{B}$. By definition, it contains  a topologically finitely generated torsion free normal closed abelian subgroup $K$  of finite index, with the quotient $G/K$ denoted by $Q$, a finite group equipped with the discrete topology. We have just seen that we have a  closed embedding $G \hookrightarrow Q \wr K$.  
\begin{Lemma} \label{abclass}For any topologically finitely generated and torsion free  profinite abelian group $A$, there is a closed embedding $A \rightarrow \hat{\bb{Z}}^n$ for some integer $n$.  
\end{Lemma} 
\begin{Proof} In Theorem 4.3.3 of  \cite{ribes}, it is shown that any torsion free topologically finitely generated abelian profinite group $A$ is of the form 
$$ \prod _p \bb{Z}_p^{r_p}
$$
where the set of integers $r_p$ are uniformly less than a fixed number $n$.  It is therefore clear that $A$ embeds in $\hat{\bb{Z}}^n \cong \prod _p \bb{Z}_p^n$.
\end{Proof}

Lemmas \ref{wreathembed} and \ref{abclass} above now give the following result. 

\begin{Proposition} For any profinite group $G$ within the class $\frak{B}$, $G$ embeds as a closed subgroup of $\Sigma _N \ltimes \hat{\bb{Z}}^N$ for some $N$.  
\end{Proposition}
\begin{Proof}  Let $G$ fit into an exact sequence of the form 
$$ 1 \rightarrow K \rightarrow G \rightarrow Q \rightarrow 1
$$
where $K$ is torsion free, abelian, and topologically finitely generated, and $Q$ is finite.  Lemma \ref{wreathembed} now shows that $G$ embeds as a closed subgroup in $Q \wr K$, and Lemma \ref{abclass} shows that $Q \wr K$ in turn embeds in $Q \wr \hat{\bb{Z}}^n$  for some $n$.  By the definition of the wreath products, it is clear that there is an embedding $Q \wr \hat{\bb{Z}}^n \hookrightarrow 
\Sigma _{q} \ltimes (\hat{\bb{Z}}^n)^{q} \cong \Sigma _q \ltimes \hat{\bb{Z}}^{nq}$, where $q = \# (Q)$.  Finally, it is also clear that the action of $\Sigma _q$ on $\hat{\bb{Z}}^{nq}$ extends over the inclusion  $\Sigma _q \hookrightarrow \Sigma _{nq}$, any permutation of the set with $nq$ elements will act on $\hat{\bb{Z}}^{nq}$ by permutation of factors.  This gives the result, with $N = nq$.  
\end{Proof}

We now use this result  to observe that members of  the family $\frak{B}$ can always be realized as Galois groups in some special ring extensions, which will turn out to be very useful in applications to $K$-theory.  Fix a ground field $k$, including all roots of unity.  For $n$ a positive integer, let $\bb{T}_n$ denote the ring $k[t_1^{\pm 1}, \ldots , t_n^{\pm 1}]$.  Similarly, let $\bb{E}_n$ denote the union
$$ \bigcup _s k[t_1^{\pm \frac{1}{s}}, \ldots , t_n^{\pm \frac{1}{s}}]
$$
where in the characteristic zero case  the union is over the partially ordered set of all positive integers (respectively the set of all positive integers prime to the characteristic of $k$ in the case of finite characteristic), with  the partial order given by $s_1 \leq s_2$ if and only if $s_1|s_2$. We now define a group action of $\hat{\bb{Z}}$ on $\bb{E}_1$ in the case where $char(k) = 0$, and of the the group $\hat{\bb{Z}}(p)$, the product of all the $q$-Sylow subgroups of $\hat{\bb{Z}}$ for $q \neq p$, in the case of positive characteristic $=p$.    In the characteristic zero case, we fix an identification $\theta : \bb{Q}/\bb{Z} \rightarrow \mu _{\infty}$, where $\mu _{\infty}$ denotes the group of all roots of unity in $k$.  In the case where $k$ has characteristic $p$, the domain of $\theta $ is $\bb{Z}_{(p)}/\bb{Z}$, where $\bb{Z}_{(p)}$ denotes the localization of the integers at $p$.  Let $\alpha \in \hat{\bb{Z}} \mbox{ (or } \hat{\bb{Z}}(p))$ denote the topological generator $1 \in \bb{Z}$. In order to specify an action of $\hat{\bb{Z}}$ or $\hat{\bb{Z}}(p)$ it suffices to specify the action of $\alpha$ on the elements $t^{\frac{1}{m}}$, which we do via the formula
$$  \alpha \cdot t^{\frac{1}{m}} = \theta({\tiny \frac{1}{m}})t^{\frac{1}{m}}
$$  This constructs an action of $\hat{\bb{Z}}$ or $\hat{\bb{Z}}(p)$ on $\bb{E}_1$. For the remainder of the paper, we will let  $\Gamma _n = 
\hat{\bb{Z}}^n$ or $\hat{\bb{Z}}(p)^n$ depending on the characteristic. Forming the tensor products $\bb{E}_n$, we obtain a $\Gamma _n$ action on $\bb{E}_n$, and it is easy to verify that $\bb{E}_n^{\Gamma_n} = \bb{T}_n$. Therefore, we have the following. 

\begin{Proposition} \label{abeliancase} The ring extension $\bb{T}_n \subseteq \bb{E}_n$ is a Galois ring extension, i.e. an infinite Galois extension in the sense of \cite{magid}, Definition 22, p.97. 
\end{Proposition} 
\begin{Proof}  We will first prove that each of the extensions $\bb{E}^{\Gamma _1}_1 \subseteq \bb{E}_1^{s \cdot \Gamma _1}$ are strongly  separable in the sense of \cite{magid}, Definition 21, p. 95.  Letting $R = \bb{E}^{\Gamma _1}_1$ and $S = \bb{E}_1^{s\cdot \Gamma _1}$, we must prove that $S$ is a projective $S \column{\otimes}{R}S$-module, where $S$  is given the $S \column{\otimes}{R} S$-module structure coming from the multiplication map  for $S$.  But this is clear, since $R = k[t^{\pm 1}]$, $S = k[t^{\pm\frac{1}{s}}]$, and therefore $S$ is a free $R$-module of rank $n$ with basis $\{1,t^{\frac{1}{s}}, \ldots , t^{\frac{s-1}{s}} \}$. The algebra $\bb{E}_1$ is the colimit of the algebras $\bb{E}_1^{s\cdot \Gamma _1}$, and so $\bb{E}_1$ satisfies the requirements of  Definition 22 of \cite{magid}.   The situation for $n > 1$ follows easily from Lemmas 4.1 and 4.2 on p. 96 of \cite{magid}.  
\end{Proof}

We also observe   that the $\Gamma _n$-action extends to an action of $\Sigma _n \ltimes \Gamma _n$ on $\bb{E}_n$.  Consider any closed subgroup $G \subseteq \Sigma _n \ltimes \Gamma _n$.  We will need a criterion to determine if the the ring extension 
$\bb{E}_n^{G} \subseteq \bb{E}$ is Galois.  We recall some terminology and a result from \cite{grothendieck}.  Let $A$ be a commutative ring equipped with the discrete topology, and let a profinite group $G$ act continuously on it with fixed point subring $B$  Let $\frak{p}$ be a prime ideal of $A$, and let $G_{\frak{p}}$ denote its stabilizer, which will be called its {\em decomposition group}. $G_{\frak{p}}$ is a closed subgroup of $G$.  Let $k({\frak{p}})$ denote the field of fractions of $A/\frak{p}$.  Similarly, let $\frak{q} = \frak{p} \cap B$, and write $k(\frak{q})$ for the field of fractions of $B/\frak{q}$.  There is an evident inclusion $k(\frak{q}) \subseteq k(\frak{p})$.    The group $G_{\frak{p}}$ acts continuously on $k(\frak{p})$, when $k(\frak{p})$ is equipped with the discrete topology.  There is therefore a homomorphism from $G_{\frak{p}}$ to the Galois  group  $Gal(k(\frak{p}), k(\frak{q}))$, which is in general a profinite group.  This homomorphism is surjective.  This is stated explicitly in \cite{grothendieck} in the case where $G$ is finite, and the result in this context follows immediately from this case by passing to inverse limits over Hausdorff finite quotients.  The kernel of this homomorphism is called the {\em inertia group} of $\frak{p}$, and is clearly a closed subgroup of $G_{\frak{p}}$.  We denote it by $I_{\frak{p}}$.  The following is proved in the case of finite $G$ in \cite{grothendieck}.  The profinite case follows directly. 

\begin{Proposition} \label{criterion}  Let $A$ be a commutative ring equipped with the discrete topology, and let a profinite group $G$ act continuously on $A$.  Then the extension $A^G \subset A$ is Galois  if and only if the inertia group $I_{\frak{p}}$ is trivial for all prime ideals $\frak{p}$ of $A$. 
\end{Proposition}

We want to consider the action of  $\Sigma _n \ltimes \Gamma _n$ on $\Bbb{E}_n$, specifically the action of $G $, where $G$ is a torsion free closed subgroup of $\Sigma _n \ltimes \Gamma _n$.   

\begin{Proposition} Let $G $ be a torsion free closed subgroup of $\Sigma _n \ltimes \Gamma _n$.  Then the ring extension $\Bbb{E}_n^G \subseteq \Bbb{E}_n$ is a Galois ring extension.  
\end{Proposition} 
\begin{Proof} We let $I_{\frak{p}}(G)$ denote the inertia group of the prime $\frak{p}$  for the group action of $G$ on $\Bbb{E}_n$.  It is readily verified that $I_{\frak{p}}(G) = I_{\frak{p}} \cap G$.  We claim that  any torsion free subgroup $G \subseteq \Sigma _n \ltimes \Gamma _n$ has a non-trivial intersection with $\Gamma _n$.  For, suppose $\gamma \in G$ is any non-trivial element.  The projection of $\gamma$ in $\Sigma _n$ is an element of finite order, say $s$.  The element $\gamma ^s$ is an element in $\Gamma _n$, and it is non-trivial due to the fact that $G$ is  torsion free.  It follows that $I_{\frak{p}}(G)  = 1$, since if there were a non-trivial element $\gamma \in I_{\frak{p}}(G)$, it would follow that $I_{\frak{p}}(G) $ would contain a non-trivial element of $\Gamma _n$, which is precluded by Proposition \ref{abeliancase}, since $I_{\frak{p}}(\Gamma _n) = I_{\frak{p}} \cap \Gamma _n$.  
\end{Proof}

\begin{Remark} {\em Note the very close analogy to the results of Auslander and Kuranishi concerning the structure of fundamental groups of flat manifolds \cite{auslander}. There it is shown that if a discrete subgroup of the group of isometries of Euclidean space $\bb{E}^n$ is torsion free, then it acts freely on $\bb{E}^n$. } 
\end{Remark}

 \end{document}